\documentclass[draft]{amsart}
\usepackage{amssymb,amscd}

\renewcommand{\P}{\mathbb P}

\renewcommand{\:}{\colon}
\renewcommand{\frak}{\mathfrak}
\renewcommand{\t}{\widetilde}

\DeclareMathOperator{\Ker}{Ker}
\DeclareMathOperator{\specan}{Specan}
\DeclareMathOperator{\ord}{ord}
\DeclareMathOperator{\supp}{Supp}
\DeclareMathOperator{\di}{div}
\DeclareMathOperator{\dlim}{dir\, lim}
\DeclareMathOperator{\LF}{LF}
\DeclareMathOperator{\diag}{Diag}
\DeclareMathOperator{\res}{res} 

\newcommand{\adeg}[1]{A_{#1}\text{-}\!\deg}
\newcommand{\aord}[1]{A_{#1}\text{-}\!\ord}
\newcommand{\wdeg}[1]{\mathbf w_{#1}\text{-}\!\deg}
\newcommand{\wpdeg}[1]{\mathbf w'_{#1}\text{-}\!\deg}
\newcommand{\word}[1]{\mathbf w_{#1}\text{-}\!\ord}
\newcommand{\wpord}[1]{\mathbf w'_{#1}\text{-}\!\ord}
\newcommand{\du}[1]{\bar A_{#1}}

\newcommand{\p}{\pi\: M \to X}
\newcommand{\X}{(X,o)}
\newcommand{\Y}{(Y,o)}

\newcommand{\Z}{\mathbb Z}
\newcommand{\Q}{\mathbb Q}
\newcommand{\R}{\mathbb R}
\newcommand{\C}{\mathbb C}
\newcommand{\N}{\mathbb N}
\newcommand{\bS}{\mathbf S}
\newcommand{\x}{\mathbf x}
\newcommand{\f}{\mathbf f}
\newcommand{\bm}{\mathbf m}
\newcommand{\w}{\mathbf w}

\newcommand{\cal}{\mathcal}
\newcommand{\A}{\mathcal A}
\newcommand{\B}{\mathcal B}
\newcommand{\E}{\mathcal E}
\newcommand{\F}{\mathcal F}
\newcommand{\M}{\mathcal M}
\newcommand{\tM}{\widetilde M}
\newcommand{\tX}{\widetilde X}
\newcommand{\tY}{\widetilde Y}
\newcommand{\tD}{\widetilde D}
\newcommand{\tF}{\widetilde F}
\newcommand{\tA}{\widetilde A}

\newcommand{\ten}{\circle*{0.25}}
\newcommand{\m}{\frak m}

 \theoremstyle{plain}
\newtheorem{thm}{Theorem}[section]
\newtheorem{cor}[thm]{Corollary}
\newtheorem{lem}[thm]{Lemma}
\newtheorem{prop}[thm]{Proposition}

 \theoremstyle{definition}
\newtheorem{defn}[thm]{Definition}

\newtheorem{ass}[thm]{Assumption}
\newtheorem{cond}[thm]{Condition}

\theoremstyle{remark}
\newtheorem{rem}[thm]{Remark}
\newtheorem{clm}{Claim}

\numberwithin{equation}{section}

\newcommand{\thmref}[1]{Theorem~\ref{#1}}
\newcommand{\lemref}[1]{Lemma~\ref{#1}}
\newcommand{\corref}[1]{Corollary~\ref{#1}}
\newcommand{\proref}[1]{Proposition~\ref{#1}}
\newcommand{\remref}[1]{Remark~\ref{#1}}
\newcommand{\clmref}[1]{Claim~\ref{#1}}
\newcommand{\defref}[1]{Definition~\ref{#1}}

\newcommand{\condref}[1]{Condition~\ref{#1}}
\newcommand{\figref}[1]{Figure~\ref{#1}}
\newcommand{\sref}[1]{Section~\ref{#1}}

\begin{document}

\title[Universal abelian covers of surface
singularities]
{Universal abelian covers of certain surface
singularities}

\author{Tomohiro Okuma}

\thanks{This research was partially supported by the
Grant-in-Aid for Young Scientists (B), 
The Ministry of Education, Culture, Sports, Science and
Technology, Japan. 
}

\address{Department of Education,
Yamagata University,
1-4-12 Kojirakawa-machi, 
Yamagata 990-8560, Japan}

\email{okuma@vc.e.yamagata-u.ac.jp}



\begin{abstract}
Every normal complex surface singularity with $\Q$-homology
 sphere link has a 
 universal abelian cover. It has been conjectured by Neumann
 and Wahl  that the universal abelian cover of a rational or
 minimally elliptic 
 singularity is a complete intersection singularity
 defined by a system of ``splice diagram equations''.
In this paper we introduce a Neumann-Wahl system,
 which is an analogue of the system of splice diagram
 equations, and prove  the following. 
 
  If $\X$ is a rational or minimally elliptic singularity,
 then its universal abelian cover $\Y$ is an equisingular
 deformation of an isolated complete intersection
 singularity $(Y_0,o)$ defined by a Neumann-Wahl system.
 Furthermore, if $G$ denotes the Galois group of the 
 covering $Y \to X$, then $G$ also acts on $Y_0$ and $X$ is
 an equisingular deformation of the quotient $Y_0/G$.
\end{abstract}

\maketitle

\section{Introduction}

Let $\X$ be a normal complex surface singularity germ.
Let $\Gamma$ and  $\Sigma$ denote the resolution graph and the link
of $\X$, respectively.
We assume that $X$ is homeomorphic to the cone over $\Sigma$.
It is known that the resolution graph and
the link of the singularity determine each other
(\cite{neumann.plumbing}).  
Assume that  the link $\Sigma$ is a $\Q$-homology sphere,
or equivalently, that 
the exceptional set of a good resolution is a tree of
rational curves.
Then $H_1(\Sigma, \Z)$ is finite.
 A morphism $\Y \to \X$ of germs of normal surface
 singularities is called a {\itshape universal abelian
 covering} if it induces an unramified Galois covering  
$Y \setminus \{o\} \to X    \setminus \{o\}$ with 
covering transformation group $H_1(\Sigma,\Z)$.
By our assumption, the universal abelian
 covering $\Y \to \X$ must exist; in fact, the link of $Y$ is
 the universal abelian cover of $\Sigma$ in the topological sense.
We are interested in the analytic properties of $Y$ and 
a way to construct  $Y$ explicitly.

In the case that $\X$ is quasihomogeneous, 
Neumann  \cite{neumann.abel} proved that the universal abelian cover
$\Y$ is a Brieskorn-Pham complete intersection, by writing down the
explicit equations  from the data of $\Sigma$
(it is known that $\Gamma$ is star-shaped in this case).
Neumann and Wahl generalized   Brieskorn-Pham complete intersections
and the way to construct them, and obtained numerous
interesting results;  
see \cite{nw-uac}, 
\cite{nw-HSL}, \cite{nw-qcusp}, \cite{nw-CIuac}.
They introduced the {\itshape splice diagram equations} 
(or {\itshape forms}) associated with a
weighted tree called a {\itshape splice diagram}
satisfying the ``semigroup condition''.
From an arbitrary resolution graph corresponding to a $\Q$-homology
sphere link, we can construct a splice diagram.
Let $\tY$ denote the singularity defined by the splice diagram
equations obtained from $\Gamma$.
They proved that $\tY$ is an isolated complete intersection surface
singularity,  and that (under ``congruence condition'') 
if the equations are chosen 
so that the {\itshape discriminant
group} $G$ ($\cong H_1(\Sigma, \Z)$) for $\Gamma$ naturally acts on
$\tY$, then the quotient $\tY/G$ is a normal surface
singularity (it is called a {\itshape splice-quotient 
singularity})  with resolution graph $\Gamma$, and the
quotient morphism is the universal abelian covering. 
Neumann and Wahl conjectured that rational singularities and
minimally elliptic singularities with $\Q$-homology sphere links are
splice-quotient singularities. 
However, it is not known whether the splice diagrams obtained
from the resolution graphs of those singularities satisfy
the semigroup condition. 

In this paper we prove that 
the universal abelian cover of  a rational or minimally elliptic
singularity is a complete intersection singularity defined
by certain  special functions. 

Let $\p$ be a good resolution with the exceptional set $A$.
Under a topological condition (\condref{c:A}), 
we can associate a collection of certain 
special polynomials and a system of weights with each node of $A$.
These polynomials are quasihomogeneous with respect to the weights.
We call the union of those collections over all nodes a
Neumann-Wahl system, 
 which is an analogue of the system of splice diagram equations. 
Though the definition of a Neumann-Wahl system is very similar to 
that of splice diagram equations, 
they are not the same (see \remref{r:diff}).
We suspect that a Neumann-Wahl system is a special type of
system of splice diagram forms.  
Our first result is the following (see \thmref{t:V}).

\begin{thm}\label{t:Vintro}
Let $(V,o)$ be a singularity defined by a Neumann-Wahl system.
Then $(V,o)$ is an isolated  complete intersection surface
 singularity.
A singularity defined by functions obtained by adding
 ``higher terms''  to the Neumann-Wahl system is an
 equisingular  deformation of $(V,o)$. 
\end{thm}

We call a singularity defined by a Neumann-Wahl system  a
Neumann-Wahl complete intersection. 
If in addition a certain analytic condition (\condref{c:B}) and a
topological condition (\condref{c:C}),
 which is stronger than \condref{c:A}, are
satisfied, then  the universal abelian
cover $Y$ is an equisingular  deformation of a Neumann-Wahl complete
intersection (\thmref{t:main-general}).  The equations of
$Y$ are constructed on the resolution space $M$.
Our equations for the deformation are automatically
equivalent with respect to a natural 
action of the discriminant group $G$, and  
the action is free on nonsingular locus.

An important point is 
that \condref{c:C} and \ref{c:B} are satisfied in case $\X$ is
a rational or minimally elliptic singularity (it can be
easily verified!). 
Thus we have the following (see \thmref{t:main} and
\proref{p:esquotient}). 

\begin{thm}\label{t:mainintro}
  If $\X$ is a rational or minimally elliptic singularity, then its
 universal abelian cover $\Y$ is an equisingular deformation of a 
 Neumann-Wahl complete intersection singularity $(Y_0,o)$.
Moreover the $\X$ is an equisingular deformation of $(Y_0/G,o)$.
\end{thm}

This paper is organized as follows.  In \sref{s:pre}, we
briefly review fundamental results on the universal abelian
covers of normal surface singularities in \cite{o.uac-rat}.  
We recall there that $\mathcal O_{Y,o}$ is
isomorphic to an $\mathcal O_{X,o}$-algebra $\A:=\bigoplus _{b \in
\B}H^0(-L^{(b)})$, where $\B$ is a group isomorphic to $G$
and $L^{(b)}$ are suitable divisors on $M$.  
In \sref{s:NWS}, we first define monomial cycles.  
The semigroup of monomial cycles is naturally
isomorphic to that of monomials; the variables are
associated with the ends of $A$.  
 We show that \condref{c:C} implies \condref{c:A}, and is
satisfied for rational or minimally elliptic singularities. 
Then we introduce the Neumann-Wahl systems and the weights.
 In \sref{s:VdefbyF}, we prove a slight generalization of
\thmref{t:Vintro}; there 
we consider a system of polynomials obtained by substituting
some power of the variables into the variables of a
Neumann-Wahl system.  
We will apply some ideas in  \cite[\S 2]{nw-CIuac} for the proof.  
In the last section, we prove \thmref{t:mainintro}.  
We construct the equations by looking at certain relations
of sections of $H^0(-L^{(b)})$'s. 
We can take a good basis of the algebra $\A$ by
\condref{c:B}, and can 
explicitly obtain the relations by \condref{c:C}.

\section{Preliminaries}\label{s:pre}

In this section, we recall some fundamental results on the universal
abelian covers of surface singularities; see
\cite{o.uac-rat} for details.  

Let $\X$ be a normal complex surface singularity germ.
We assume that the link $\Sigma$ of $\X$ is a $\Q$-homology sphere
and that $X$ is homeomorphic to a cone over  $\Sigma$.
Then $X$ has a unique universal abelian cover. 
Let $\p$ be a  resolution of the singularity, and let
$A=\bigcup _iA_i$ be the decomposition of the exceptional set
$A=\pi^{-1}(o)$ into irreducible components.
Assume that $\pi$ is a good resolution, i.e., 
 $A$ is a divisor having only simple normal crossings.
Then the condition that $\Sigma$ is a $\Q$-homology sphere is
equivalent to that $H^1(\mathcal O_A)=0$, i.e., $A$ is a tree of
nonsingular rational curves.
We call a divisor supported in $A$ a cycle.
Let  $A_{\Z}$ denote the group of cycles.
An element of $A_{\Q}:=A_{\Z}\otimes \Q$   is called a $\Q$-cycle.
Let $\du i \in A_{\Q}$ denote the dual cycle of $A_i$, i.e., 
the $\Q$-cycle
satisfying  $\du i\cdot A_j=-\delta _{ij}$, 
where $\delta _{ij}$ denotes the Kronecker delta.
We denote by $\du{\Z}$ the subgroup of $A_{\Q}$ generated by
$\du i$'s. 
Recall that the first homology group $H_1(\Sigma, \Z)$,  the  Galois
group of the universal abelian covering of $X$,  is isomorphic
to the group $\du {\Z}/A_{\Z}$ called the discriminant group. 
The order of the group is $|\det(A_i\cdot A_j)|$. 
Let $D$ be a $\Q$-divisor on $M$.
 We denote by $\nu(D)$ the $\Q$-cycle
satisfying $(\nu(D)-D)\cdot A_i=0$ for all $A_i$.
We say that $D$ is $\pi$-anti-nef if $-D$ is $\pi$-nef,
i.e., $D\cdot A_i \le  0$ for all $A_i$.
Let  $ \F (D)$ denote the set of $\pi$-anti-nef $\Q$-divisors $F$
 satisfying  $F-D\in A_{\Z}$. 
Note that $\F(D)$ has the minimum with respect to ``$\ge$''.

We take effective $\Q$-cycles $E_1, \ldots , E_s$ such that
if $\E_i$   
denotes the cyclic subgroup of $\du {\Z}/A_{\Z}$ generated by $E_i$,
then $\du {\Z}/A_{\Z}=\E_1 \oplus \cdots \oplus \E_s$.
Let $r_i$ be the order of $\E_i$. 
Then for $1 \le i \le s$ there exists a divisor 
$L_i$ and a function $f_i$ on a suitable neighborhood 
of $A$  such that $r_iL_i- r_iE_i=\di (f_i)$. 
For any $b=(b_1,\dots,b_s)\in \Z^s$, we define a divisor
$L^{(b)}$ by  
$$ 
L^{(b)}=\sum_{j=1}^s b_jL_j-\left[\sum_{j=1}^sb_jE_j\right].
$$
Let  $\bar b_i$ denote the smallest
nonnegative integer such that $r_i \mid b_i -\bar b_i$.
Let $\bar b= (\bar b_1, \dots , \bar b_s)$ and 
 $\B=\{\bar b| b \in \Z^s\}$; the set $\B$ is identified with the
 discriminant group $\du
 {\Z}/A_{\Z}$. 
We define  a  sheaf $\bar \A$ of $\mathcal O_M$-modules by 
$$
\bar \A=\bigoplus _{b \in \B}\mathcal O_M(-L^{(b)}).
$$
The $\mathcal O_M$-algebra structure of $\bar \A $ is given by the composite
$$
\mathcal O_M(-L^{(b)}) \otimes \mathcal O_M(-L^{(b')})  \to \mathcal O_M(-L^{(b)}-L^{(b')})
\subset  \mathcal O_M(-L^{(b+b')})
$$
and the isomorphism
$$ 
 \mathcal O_M(-L^{(b+b')}) \to \mathcal O_M(-L^{(\overline{b+b'})})
$$
 given by multiplying by $\prod
_{b_i+b_i'\neq\overline{b_i+b'_i}}f_i^{-1}$. 
Then the natural projection 
$$
 Y:=\specan_X \pi_*\bar\A \to X
$$
is the universal abelian covering 
(see \cite[Theorem 3.4]{o.uac-rat}). 
The local ring $\mathcal O_{Y,o}$ of the singularity $\Y$ is 
isomorphic to
$$
\A:=(\pi_*\bar\A )_o=\bigoplus _{b \in \B}H^0(-L^{(b)}),
$$
where 
$H^0(-L^{(b)})=\dlim _U H^0(U,\mathcal O_M(-L^{(b)}))$,
$U$ varies over all open neighborhoods of $A$.
We write $\A_b=H^0(-L^{(b)})$.

Let $h \in \A_b$. 
If  $\di (h)-L^{(b)}-D$ has no component of $A$ for some
cycle $D\in A_{\Z}$, 
then we write $(h)_A=\nu(L^{(b)})+D$.

\begin{lem}\label{l:product}
 Let $\sigma_i \in \A_{b^i}$, $i=1,2$, and 
let $\sigma_1\cdot \sigma_2\in \A_{\overline{b^1+b^2}}$ be
 the product of $\sigma_1$ and $\sigma_2$   in the algebra $\A$. 
Suppose that a divisor $L \in \F(L^{( \overline{b^1+b^2})})$
 satisfies 
 $\nu(L)=(\sigma_1)_A+(\sigma_2)_A$. 
Then $\sigma_1\cdot \sigma_2$
is a section of  $H^0(-L) $
and $(\sigma_1\cdot \sigma_2)_A=\nu(L)$. 
\end{lem}
\begin{proof}
 It follows from the definition of the algebra structure of
 $\A$. 
\end{proof}

 Let $D$ be a reduced and connected cycle.
A component $A_i$ of $D$ is called an {\itshape end} of $D$
if $(D-A_i)\cdot A_i\le 1$.
We denote by $\E(D)$ the set of ends of $D$.

 Assume that $\E(A)=\{A_1, \dots ,A_m\}$.
For any $1 \le i \le m$, there uniquely exists a divisor
$L^i$ such that $\nu
(L^i)=\du i$ and $L^i \in \F(L^{(b)})$ for some $b \in \B$.
If $\X$ is rational, then there exists $y_i \in H^0(-L^i)$
such that $(y_i)_A=\nu(L^i)$ (cf. \lemref{l:satisfyB}).
The next theorem follows from  \cite[Theorem 7.5]{o.uac-rat} and its
proof.

\begin{thm}\label{t:repOY}
Assume that $\X$  is rational.
 Let $y_i \in H^0(-L^i)$, $1 \le i \le m$, be  as above.
 Let  $S=\C\{x_1, \dots ,x_m\}$ be the convergent power series
ring. 
We define a homomorphism $\psi \: S \to \A=\mathcal O_{Y,o}$ of
 $\C$-algebras by $\psi(x_i)=y_i$. 
Then $\psi$ is  surjective.
\end{thm}

In the last section, we give a proof of \thmref{t:repOY},
which is different from that in \cite{o.uac-rat}.
In fact it is shown  that the
assertion also holds true for minimally elliptic singularities.

\section{Neumann-Wahl systems}\label{s:NWS}

In this section we will introduce a Neumann-Wahl system
associated with the exceptional set $A$.
It is a set of certain polynomials, 
and an analogue of the system of splice diagram
equations in Neumann and Wahl's work (\cite{nw-uac}, \cite{nw-HSL},
\cite{nw-CIuac}).

We use the notation of the preceding section, and keep
the assumption that $H^1(\mathcal O_A)=0$.
First assume that the set $\E(A)$ consists of the
components $A_1, \dots ,A_m$. 
Let $\C[x_1, \ldots ,x_m]$ be the polynomial ring.
 \thmref{t:repOY} suggests us the following

\begin{defn}
 Let $D=\sum a_i\du i \in A_{\Q}$, $a_i \ge 0$. 
If $a_i=0$ for all $A_i \notin \E(A)$,
then we call $D$ a {\itshape $\Q$-monomial cycle};
if in addition $a_i \in \Z$ for all $i$, we call $D$ a
 {\itshape monomial cycle}. 
 For any monomial cycle $D=\sum _{i=1}^ma_i\du i$, we associate a
monomial 
$$
 x(D):=\prod _{i=1}^mx_i^{a_i} \in \C[x_1, \ldots ,x_m].
$$
The $x$ induces an isomorphism between the semigroup of
 monomial cycles and that of monomials of $x_1,
 \ldots, x_m$. 
Formally, we may also consider the $\Q$-monomial $x(D)$ for a
 $\Q$-monomial  cycle $D$. 
\end{defn}

\begin{defn}
For any $F=\sum a_kA_k \in A_{\Q}$, we write $m_{A_k}(F)=a_k$.
 For any component $A_j$, we define the {\itshape
 $A_{j}$-weight} of $x_i$ to be  $m_{A_{j}}(\du i)$.
Then the {\itshape $A_{j}$-degree} of a monomial $x(D)$ is
 defined to be $m_{A_{j}}(D)$. 
\end{defn}

A connected component of 
$A-A_i$ is called a {\itshape branch} of $A_i$.
A component $A_i$ is called a {\itshape node} if $(A-A_i)\cdot
A_i\ge 3$. 
We consider the following conditions concerning the weighted
dual graph of $A$.

\begin{cond}\label{c:A}
For  any branch $C$ of any node $A_i$,
 there exists a monomial cycle $D$
 such that $D-\du i$ is an effective integral cycle
 supported on $C$;  
in this case, we say that  $D$ (or monomial
$x(D)$) belongs to the branch $C$.
\end{cond}

Note that in general there may exist more than one monomials
belonging to a branch.

\begin{cond}\label{c:C}
$A$ is star-shaped, or for any branch $C$ of any component
 $A_i\notin \E(A)$,  the fundamental cycle $Z_C$ supported
 on $C$ satisfies $Z_C\cdot  A_i=1$.
\end{cond}

\begin{lem}\label{l:existmonomials}
\condref{c:C} implies \condref{c:A},
and is satisfied in the following cases:
\begin{enumerate}
 \item $\X$ is a rational singularity;
 \item $\X$ is a minimally elliptic singularity, and
 the minimally elliptic cycle is supported on $A$
(this condition is satisfied on the minimal good  resolution).
\end{enumerate}
\end{lem}
\begin{proof}
Assume that the condition (1) or (2) is satisfied.
From  basic results on the computation sequences for the
 fundamental cycle
 (\cite{la.rat},  \cite{la.me}), we obtain
\condref{c:C}.
Suppose that $C_{1}$ is a branch of a node $A_{i_1}$. 
Then $D_1:=\du {i_1}+Z_{C_1}$ is $\pi$-anti-nef and $D_1 \cdot A_j=0$ 
for every $A_j \le A-C_1$. 
If $D_1 \cdot A_{i_2}<0$ for $A_{i_2} \notin \E(A)$, then take a 
 branch $C_2$ of $A_{i_2}$, not
 containing $A_{i_1}$, and put $D_2:=D_1+Z_{C_2}$.
In this manner we obtain  a finite sequence $\{D_1, \ldots ,D_n\}$ 
of  $\pi$-anti-nef cycles, which ends with a monomial cycle
 belonging to $C_1$. Thus \condref{c:A} is satisfied.
These arguments also show that \condref{c:A} holds in case
 $A$ is star-shaped.
\end{proof}

\begin{defn}\label{d:CSAM}
Assume that  \condref{c:A} is satisfied, and that
  $A_1, \ldots ,A_s$ are all of the nodes of $A$.
Let $C_1, \dots ,C_p$ be the branches of $A_{1}$.

(1) A monomial $x(D)$ is 
called an {\itshape admissible monomial} at the node $A_{1}$
if it belongs to one of the branches of $A_1$. 
 A set of monomials $\{x(D_1), \dots ,x(D_p)\}$ is called a
 {\itshape complete system} of admissible monomials 
at $A_{1}$ if $D_i$ belongs to $C_i$ for 
 $i=1, \ldots ,p$.

(2)  Let $\{m_1, \ldots ,m_p\}$  be any complete system of
 admissible  monomials  at $A_{1}$.
Let $F=(c_{ij})$, $c_{ij} \in \C$, 
be a $((p-2) \times  p)$-matrix such  
 that every maximal minor of it has rank $p-2$.
We define polynomials $f_1, \ldots ,f_{p-2}$ by 
$$
\begin{pmatrix}
 f_1\\ \vdots \\ f_{p-2}
\end{pmatrix}
= F
\begin{pmatrix}
 m_1\\ \vdots \\ m_p
\end{pmatrix}.
$$
 We call each $f_i$ an {\itshape admissible form} at $A_1$ and
the set $\{f_1 , \dots ,f_{p-2}\}$ a {\itshape Neumann-Wahl
 system}  at  $A_{1}$. 

(3) Let $\F_i$ denote a Neumann-Wahl system  at a node  $A_{i}$.
Then we call the set 
$ \bigcup _{i=1}^s\F_i$ a {\itshape Neumann-Wahl system} associated
 with $A$; it is an empty set in case $A$ has no nodes. 
\end{defn}

\begin{rem}\label{r:normal}
The matrix $F$ above  can be reduced to the following matrix by
 row operations:
$$
\begin{pmatrix}
  1 & 0 & \dots & 0 & a_1 & b_1 \\
0 & 1 & \dots & 0 & a_2 & b_2 \\
\vdots & \vdots & \ddots  & \vdots  & \vdots  & \vdots  \\
 0 & 0 & \dots  & 1 & a_{p-2} & b_{p-2}
 \end{pmatrix},
$$
where $a_ib_j-a_jb_i\neq 0$ for $i\neq j$, and all $a_i$ and $b_i$
 are nonzero.
\end{rem}

The admissible forms at a node $A_i$ are 
quasihomogeneous polynomials with respect to the
$A_i$-weight.
The following lemma is needed in the next section.

\begin{lem}\label{l:higher}
 Let $D$ be a $\pi$-anti-nef $\Q$-cycle such that
 $m_{A_i}(D)>m_{A_i}(\du i)$ for some $i$. 
Then for any component $A_j$, we have $m_{A_j}(D)>m_{A_j}(\du i)$.
\end{lem}
\begin{proof}
First we will show that $D \ge \du i$.
There exist effective $\Q$-cycles $D_1$ and $D_2$ such that 
$D-\du i=D_1-D_2$ and that $D_1$ and $D_2$ have no common
 components.
If $D_2>0$, then
$$
m_{A_i}(D_2)=-\du i\cdot D_2=-D\cdot D_2+D_1\cdot
 D_2-D_2\cdot D_2>0.
$$
It contradicts the assumption of the lemma. Hence $D-\du
 i\ge 0$. 

 Let $a=m_{A_i}(D-\du i)$ and $C$ a branch of  $A_i$.
If a component $A_k$ of $C$ intersects $A_i$, then $(\du
 i+aA_i)\cdot  A_k=a>0$. 
Since $D$ is $\pi$-anti-nef, there exists a $\Q$-cycle $C'$
 supported on  $C$ such that $C'\cdot A_l \le 0$ for all
 $A_l \le C$ and that  
$D \ge \du i+aA_i+C'$.
\end{proof}

\begin{rem}\label{r:diff}
The  definition of  a Neumann-Wahl system is very similar to that of
a system of splicing diagram equations in the Neumann and Wahl's work.
However those are not the same.
Let us consider a weighted graph $\Gamma$ represented as in
 \figref{f:Gamma}; 
the vertex  \rule{2mm}{2mm} has weight $-4$ and other vertices
$\bullet$ have  weight $-2$.
\begin{figure}[htbp]
  \begin{center}
\setlength{\unitlength}{0.5cm}
    \begin{picture}(11,2)(0,0)
\multiput(0,0)(0,2){2}{\ten}
\multiput(5,1)(2,0){2}{\ten}
\put(2.8,0.8){\rule{2mm}{2mm}}
\multiput(9,0.5)(0,1){2}{\ten}
\multiput(11,0)(0,2){2}{\ten}
\put(3,1){\line(1,0){4}}
\put(0,0){\line(3,1){3}}
\put(0,2){\line(3,-1){3}}
\put(7,1){\line(4,1){4}}
\put(7,1){\line(4,-1){4}}
\end{picture}
\caption{\label{f:Gamma}}
  \end{center}
\end{figure}
The graph $\Gamma$ is realized as the
 resolution graph of an elliptic 
 singularity with $\Q$-homology sphere link.
The splice diagram $\Delta$ associated with $\Gamma$ is
 represented  as  in \figref{f:Delta}.
\begin{figure}[htbp]
  \begin{center}
\setlength{\unitlength}{0.5cm}
    \begin{picture}(10,2)(0,0)
\multiput(0,0)(0,2){2}{\ten}
\multiput(3,1)(4,0){2}{\ten}
\multiput(10,0)(0,2){2}{\ten}
\put(3,1){\line(1,0){4}}
\put(0,0){\line(3,1){3}}
\put(0,2){\line(3,-1){3}}
\put(7,1){\line(3,1){3}}
\put(7,1){\line(3,-1){3}}
\put(2.2,1.5){$2$} 
\put(2.2,0){$2$} 
\put(7.5,1.5){$3$} 
\put(7.5,0){$3$} 
\put(3.5,1.2){$3$} 
\put(6,1.2){$20$} 
\end{picture}
\caption{\label{f:Delta}}
  \end{center}
\end{figure}
We see that $\Delta$ satisfies semigroup condition.
However,  $\Gamma$ does not  satisfy \condref{c:A}.
Therefore we cannot define Neumann-Wahl systems in this
 case, though the  splice diagram 
 equations are  defined. 
That might indicate a Neumann-Wahl system is a special type
 of 
system of  splice diagram equations.
\end{rem}

\section{Varieties defined by Neumann-Wahl systems}\label{s:VdefbyF}

In this section we prove that any Neumann-Wahl system
 defines a complete intersection
surface with an isolated singularity at the origin; 
so we will call such a singularity a {\itshape Neumann-Wahl complete
 intersection singularity}.
In fact, we will prove the assertion for a slight generalization of
a Neumann-Wahl system.
We note that some of   methods in this section  are  discussed  in
\cite{nw-CIuac}. 

We use the notation of the preceding section.
Assume that \condref{c:A} is satisfied and that $A$ has at
least one node.
As in the preceding section, we associate  ends of $A$ with 
the variables $x_1, \ldots
,x_m$, where   $m=\# \E(A)$. 
Suppose that $A_1, \ldots ,A_s$ are all of the nodes of $A$.
Let $d_i$ denote  the number of branches of a node $A_i$.
Let $\M_i=\{m_{i1}, \ldots ,m_{id_i}\}$ denote a complete system of
admissible monomials and $\F_i=\{f_{i1}, \ldots ,f_{id_i-2}\}$
a Neumann-Wahl system at a node $A_i$, where each $f_{ij}$
is a linear form of monomials of $\M_i$.
By counting the numbers of the ends, the nodes and the edges
around the nodes of the dual graph of $A$, we   see that
$$
\sum _{i=1}^s(d_i-2)=m-2.
$$

Let $C_1, \ldots ,C_{d_1}$ denote the branches of $A_1$.
Without loss of generality, we may assume the following.
\begin{enumerate}
 \item For $1\le j \le d_1-1$, $C_j$ is a chain of curves; 
in this case, $A_1$ is
 an end of the minimal reduced connected cycle containing
       all nodes of $A$.
 \item For  $1\le j \le d_1-1$, the variable $x_j$
 corresponds to the end of $C_j$.
 \item For  $2\le i \le s$, the admissible monomial $m_{id_i}$
       belongs to the  branch of $A_i$ containing $A_1$.
 \item For $1\le i \le s$, the admissible  forms of $\F_i$ are
       given by the 
       matrix as  in \remref{r:normal};  we write
 \begin{align*}
 f_{ij} 
&= m_{ij}+a_{ij}m_{id_i-1}+b_{ij}m_{id_i},  & &1\le i \le s, \;
1\le j \le d_i-2, \\
 m_{1j} & =  x_j^{\alpha_j}, & &1\le j \le d_1-1. 
\end{align*}
\end{enumerate}

Now we slightly modify the admissible forms.
This modification is needed for the induction step of the
proof of the main theorem.
Let $\N$ denote the set of positive integers. 
A vector $v \in \N^m$ is said to be {\itshape primitive} if
$v$ cannot be written as $v=c v'$ 
with $v' \in \N^m$ and $c \in \N$, $c>1$.
Fix an arbitrary vector $\delta=(\delta_1, \ldots
,\delta_m) \in \N^m$. 
For each node $A_i$, let $e_i$ denote the positive integer
such that 
$$
\w_i:=(\adeg i(x_1)\cdot e_i/\delta_1, \ldots ,
\adeg i(x_m)\cdot e_i/\delta_m) \in \N^m
$$
is primitive.
 Let  $S=\C\{x_1, \dots ,x_m\}$ be the convergent power series
ring. 

\begin{defn}\label{d:wdeg}
Let $\w =(w_1, \ldots ,w_m)\in \N^m$. 
We define the {\itshape $\w$-degree} of a $\Q$-monomial $m=\prod
 _{k=1}^mx_k^{a_k}$ to be  $\wdeg{}(m)=\sum a_kw_k$.
Let $f=\sum_{k\ge 1} f_k \in S$, where $f_1\neq 0$
 and each $f_k$ is a quasihomogeneous polynomial with respect to 
$\w$ such that  $\wdeg {}(f_k)<\wdeg{}(f_{k+1})$.
We call $f_1$ the {\itshape leading form} of $f$, and 
 denote it by $\LF_{\w}(f)$.
Then $f-\LF_{\w}(f)$ is called  the
 {\itshape higher term} of $f$. 
We define {\itshape $\w$-order} of $f$ to be 
$\mathbf w\text{-}\!\ord (f)=\mathbf w\text{-}\!\deg
(\LF_{\w}(f))$. We set $\mathbf w\text{-}\!\ord (0)=\infty$.
\end{defn}

We write $\x_k=x_k^{\delta_k}$ and $\f=f(\x_1, \ldots ,\x_m)$ for
$f=f(x_1, \ldots ,x_m) \in S$.
We call $\f$ the {\itshape $\delta$-lifting} of $f$.
Any monomial can be thought as the $\delta$-lifting of a
$\Q$-monomial. 
A polynomial $f$ is quasihomogeneous with respect to  $A_i$-weight
if and only if so is $\f$ with respect to the weight $\w_i$; in
fact, for a $\Q$-monomial $m$, we have 
\begin{equation}\label{eq:deg}
 \wdeg i(\bm)=\adeg
i(m)\cdot e_i.
\end{equation}

For each $\f_{ij}$, we take a convergent power series
$f_{ij}^+ \in  S$
satisfying 
$$
 \wdeg i(\f_{ij})<\word i(f_{ij}^+).
$$
For each $t \in \C$, we set  
$$
 \F_t=\{\f _{ij}+tf_{ij}^+ | 1 \le i \le s, \; 1\le j \le d_i-2\}.
$$
We note that if $m_{ij}=x(D_{ij})$, then for $i \ge 2$ and
$1 \le j \le d_i-1$, 
$$
 m_{A_1}(D_{id_i})>m_{A_1}(D_{ij})=m_{A_1}(\du i).
$$
If the $\delta$-lifting of a $\Q$-monomial $x(D)$ is
contained in $f_{ij}^+$, then $m_{A_i}(D)>m_{A_i}(\du i)$.
By \lemref{l:higher}, $m_{A_1}(D)>m_{A_1}(\du i)$.
Thus we obtain  the following:
\begin{align*}
 \LF_{\w_1}(\f_{1j}+tf_{1j}^+)&=\f_{1j}, \\
 \LF_{\w_1}(\f_{ij}+tf_{ij}^+)& =\bm_{ij}+a_{ij}\bm_{id_i-1}, 
 \quad \text{for $i \ge 2$}.
\end{align*}
Therefore the set  
$$
 \LF_{\w_1}\F:=\{\LF_{\w_1}(f)|f \in \F_t\}
$$
 is independent of $t \in \C$; it can be shown that
 $\LF_{\w_i}\F$ is also 
 independent of $t$ for $2 \le i \le s$.

\begin{defn}[{Wahl \cite{wahl.es}, cf. \cite[V]{la.simul}}]
 Let $\omega \: \tX \to T$ be a deformation of a normal surface
 singularity  $\tX_o=\omega^{-1}(o)$, $o \in T$. 
Suppose that each fiber $\tX_t$ has only one singular point
 and that  
there exists a simultaneous resolution  $\bar \omega\: \tM \to \tX$
 with the exceptional set $\tA$.
If the restriction $(\omega \circ \bar \omega)|_{\tA}$ is
 a locally  trivial deformation of the exceptional divisor
 of $\tM_o$, then we call
 $\omega \circ \bar \omega$ (resp. $\omega$) an {\itshape
 equisingular  deformation}.
$\bar \omega$ is called a {\itshape weak simultaneous resolution} of
 $\omega$. 
\end{defn}

For a subset $B$ of any commutative ring, 
let $I(B)$ denote the ideal generated by the elements of $B$.

Let $(V_t, o) \subset (\C^m,o)$ denote the singularity
defined by the ideal $I(\F_t) \subset S$. 
The main result of this section is the following.

\begin{thm}[cf. {\cite[Theorem 2.6]{nw-CIuac}}]\label{t:V}
The singularity $(V_t ,o)$ is an isolated  complete
 intersection surface  singularity  for each $t \in \C$.
Furthermore, the family $\{V_t | t \in \C\}$ is an equisingular
 deformation.
\end{thm}

First we show that every $V_t$  is a complete intersection.

\begin{lem}[cf. {\cite[Theorem 3.1]{nw-CIuac}}]\label{l:curveC}
For any variable $x_k$, let $A_{i_k}$ denote the node
 nearest to the end corresponding to $x_k$. 
 Let $C \subset \C^{m}$ be the affine variety defined by 
the ideal $I(\LF_{\w_{i_k}}\F\cup \{x_k\})$. 
\begin{enumerate}
 \item If $x_j=0$ ($j\neq k$) at $p \in C$, then $p$ is the origin.
 \item $C$ is a complete intersection curve and smooth
       except for the origin. 
\end{enumerate}
\end{lem}
\begin{proof}
Without loss of generality, we may assume that $k=1$.
Let 
$$
c_{j}=\begin{cases}
       b_{11}/a_{11} & \text{if $j=d_1-1$}, \\
       b_{1j}-a_{1j}b_{11}/a_{11} & \text{if $2 \le j\le d_1-2$}.
      \end{cases}
$$
Then $c_j\neq 0$ and  $C$ is a subvariety of the affine
 space $\C^{m-1}$ with  coordinates  $x_2, \dots, x_m$ 
 defined by the equations
\begin{equation}\label{eq:C}
 \begin{split}
   \x_j^{\alpha_j}+c_j\bm_{1d_1}=0, \quad  & 2 \le j \le d_1-1, \\
  \bm_{ij}+a_{ij}\bm_{id_i-1}=0, \quad    & 2\le i \le s,
  1\le j \le 
 d_i-2.
 \end{split}
\end{equation}
If a monomial appearing in  \eqref{eq:C}
 vanishes at $p \in C$, 
 then so does every monomial in the equations at the same node.
On the other hand, for each $2 \le i \le m$, some power of
 $x_i$ appears  in  \eqref{eq:C}  
because each end of $A$ is a unique end of a branch of the
 nearest node. 
Thus if a variable $x_i$ ($i \ge 2$) vanishes on $C$, then
 so do all monomials appearing in \eqref{eq:C}, since $A$ is
 a connected tree of curves.  
Hence we obtain (1).

Since  $\#\LF_{\w_1}\F=m-2$,   it follows from (1) that
 $\{x_1, x_i\} \cup  \LF_{\w_1}\F$, where $i\neq 1$, is a 
 regular sequence.
Hence $C$ is one-dimensional and complete intersection.
The argument above also shows that an ideal $I(\{x_1, x_i-1\} \cup
 \LF_{\w_1}\F)$ 
 defines  a nonsingular zero-dimensional variety.
Since $C $ is defined by
 quasihomogeneous polynomials, $C$ is smooth except for  the
 origin. 
\end{proof}

\begin{cor}[cf. {\cite[Corollary 3.4]{nw-CIuac}}]\label{c:CI}
 $V_t$ is a complete intersection surface singularity, and
 the support  of $V_t \cap \{x_j=x_k=0\}$, $j\neq k$, 
is the origin.
\end{cor}
\begin{proof}
By \lemref{l:curveC}, $\LF_{\w_{i_k}}\F\cup \{x_j,x_k\}$ is
 a regular  sequence. 
Hence so are $\F_t$ and $\F_t\cup\{x_j,x_k\}$.
\end{proof}

\begin{defn}
 We define the {\itshape weighted dual graph} of a normal
 surface singularity to be
 that of the 
 exceptional set of the minimal good resolution of the singularity.
\end{defn}

Let $(W,o)$ be a germ of a normal surface singularity and
$W' \to W$ the minimal resolution. 
The {\itshape canonical cycle} on $W'$ is a
$\Q$-cycle which is 
numerically equivalent to the canonical divisor $K_{W'}$.
The self-intersection number of the canonical cycle is an
invariant of the singularity and 
determined by the weighted dual graph; we denote it by $K^2(W)$.
Let $\omega \: \tX \to T \subset \C$ be  
a deformation of surface singularities.
The invariance of $K^2(\tX_t)$ implies the existence of the
simultaneous canonical model (or simultaneous RDP
resolution) of $\omega$;  first the  Gorenstein case was proved by
Laufer  (\cite[Theorem 4.3]{la.simul}),  
and the general case by Ishii (\cite[Corollary 1.10]{is.simul}).
By  \cite{bri.simul}, the singularities of the simultaneous
canonical 
model are simultaneously resolved after a suitable finite
 base change. 
Thus  we obtain the following theorem by the  arguments of
\cite[VI]{la.simul}. 

\begin{thm}[Laufer, Ishii]\label{t:la}
  Let $\omega \: \tX \to T\subset \C$ be a deformation of a normal
 surface  singularity. If the weighted dual graphs of
 $\tX_t$, $t\in T$,  are the same, then $\omega$ is  an
 equisingular deformation, and it 
 admits a simultaneous resolution such that each fiber is
 the minimal  good resolution.
\end{thm}

For a divisor $D$, we denote by $D_{red}$ the reduced
divisor with $\supp (D_{red})=\supp (D)$. 
For the induction step of the proof of the main theorem,
 we need the following.

\begin{lem}\label{l:es-div}
  Let $\omega \: \tX \to \C$ be an equisingular
 deformation of a germ of a normal surface
 singularity  $\tX_0=\omega^{-1}(0)$.
Let $\tD$ be a reduced divisor on $\tX$, which contains the
 singular locus of $\tX$.
Suppose that $\omega |_{\tD}$ is a locally
 trivial  deformation of a reduced divisor $\tD_0:=\t D|_{\tX_0}$ on
 $\tX_0$. 
Then there exists a simultaneous resolution  $\bar \omega\:
 \tM \to \tX$  with the exceptional set $\tA$ such that  
$(\omega \circ \bar \omega)|_{(\bar \omega^*\tD)_{red}}$ is
 a locally  trivial deformation 
(hence so is $(\omega \circ \bar \omega)|_{\tA}$).
\end{lem}
\begin{proof}
 There exists a weak simultaneous resolution $\bar \omega\:\tM
 \to \tX$ with  the exceptional set $\tA$ such that each
 $\tA_t:=\tA |_{\tM_t}$ has only simple normal 
 crossing and $(\omega \circ \bar \omega)|_{\tA}$ is a
 locally trivial  deformation. 
Let $\tF$ denote the strict transform of $\tD$ on $\tM$.
Let $\tA^1$ (resp. $\tF^1$) be a divisor on $\tM$, which
 is the total space of the deformation of an irreducible
 component of $\tA_0$ (resp. $\tF_0:=\tF |_{ \tM_0}$).
Then the intersection number $\tA^1_t\cdot \tF_t^1$ is constant.
We may assume $\# (\tA^1_t\cap \tF_t^1) \le 1$.
If $\tA^1_t\cdot \tF_t^1\ge 2$, then take the blowing up of
 $\tM$ along  the curve $\tA^1\cap \tF^1$. 
By taking blowing ups successively in a similar way, we obtain a
 simultaneous resolution $\bar \omega'\: \tM' \to \tX$ such that 
each divisor $((\bar  \omega'|_{\tM'_t})^*\tD_t)_{red}$ has
 only normal  crossings 
and that the  weighted dual graph of the divisor is
 independent of $t \in \C$. 
\end{proof}

\begin{lem}\label{l:G}
 Let $\omega\: \tX \to \C$ be as in \lemref{l:es-div}.
Suppose that $\tX$ is embedded in an open subset of
 $\C^n\times \C$ such  that the singular locus of $\tX$ is
 $\{o\}\times \C$, and that $\omega$  is the composite of
 this embedding and the projection $\C^n\times \C \to \C$.  
Let $G$ be a finite subgroup of the unitary group $\mathrm
 U(n)\subset \mathrm {GL}(\C^n)$.  
Then $G$ acts on $\C^n\times \C$ by $g\cdot(z,t)=(g\cdot z,
 t)$, $g \in G$.
Assume the action  induces an action on
 $\tX$ which is free on  $\tX \setminus \{o\}\times \C$.
Then the morphism $\mathcal Omega\: \tX /G \to \C$ obtained from $\omega$
 is an  equisingular deformation of $(\tX_t/G,o)$, $t \in \C$. 
\end{lem}
\begin{proof}
Let $\bS_c \subset \C^n$, $c>0$, denote the $(2n-1)$-sphere
 of radius $c$. 
Let $t_0 \in \C$ be an arbitrary point.
 Then there exist an open neighborhood $U$ of $t_0$ and  a
 positive number  $\epsilon \in \R$ such that  
$\Sigma_t:=\bS_{\epsilon}\cap \tX_t \subset \C^n$ is the
 link of $\tX_t$  for every $t \in U$ and the family 
$\{\Sigma_t | t \in U\}$ is topologically trivial. 
By the assumption, $G$ acts on $\{\Sigma_t | t \in U\}$ freely.
Thus  we obtain a family $\{\Sigma_t/{G} | t \in U\}$ which is 
topologically trivial.
Recall that the weighted dual graph of a surface singularity is
determined by its link (Neumann \cite{neumann.plumbing}).
By \thmref{t:la}, we obtain the assertion.
\end{proof}

We mention the weighted blowing up which is needed in the
proof of the theorem.  
Let $\w=(w_1, \ldots ,w_m) \in \N^m$ be a primitive vector, and 
let $\beta\: Z \to \C^m$ be the weighted blowing up with
respect to the  weight $\w$.
It is a projective morphism  inducing an isomorphism 
$Z \setminus \beta^{-1}(o) \to \C^m \setminus  \{o\}$,   and
$\beta^{-1}(o)=\P(\w)$, the weighted projective space of
type $\w$. 
The variety $Z$ is covered by affine varieties $Z_1, \dots ,Z_m$; 
each  $Z_i$ is a quotient of $W_i=\C^m$ by a cyclic group $\cal C_i$
of
order $w_i$  determined by the weight $\w$.
Let $\{x_1, \ldots ,x_m\}$ and $\{z_1, \ldots ,z_m\}$
be the coordinates of $\C^m$ and  $W_1$, respectively. 
The action of the group $\cal C_1$ on $W_1$ is given by the diagonal
matrix 
$$
\diag [e(-1/w_1), e(w_2/w_1), \, \ldots \, ,e(w_m/w_1)],
$$
where  $e(q)=\exp(2\pi \sqrt{-1}q)$. Since $\w$ is
primitive, $\cal C_1$ is trivial or the fixed
locus is a proper subvariety of the hyperplane $\{x_1=0\}$
which is the exceptional locus.
The morphism $W_1  \to \C^m$, 
which is the composite of the quotient morphism $W_1 \to Z_1$ 
and $\beta\:Z_1 \to \C^m$,
 is given  by 
$$
x_1=z_1^{w_1}, \quad x_i=z_1^{w_i}z_i \quad (i=2, \ldots ,m).
$$

\begin{proof}[Proof of \thmref{t:V}]
We prove the theorem by induction on the number of nodes $s$ of $A$.
We have to show the isolated singularity of each $V_t$ and the
 equisingularity of the family $\{V_t | t \in \C\}$.

First assume that $s=1$.
Then  $V_0$ is the so-called Brieskorn-Pham complete intersection
 singularity; it is known that $V_0$ has an isolated
 singularity 
(it is  also easily  checked  by using the   Jacobian criterion).
We fix $t \in \C$. 
Since $\LF_{\w_1}\F$ is a regular sequence,
it follows from the theory of filtered rings that there exists an
 equisingular deformation of $V_0$  
with  general  fiber $V_t$  (cf. \cite[\S 6]{tki-w},
 \cite{wahl.defqh}). 
Hence $V_t$ is an isolated complete intersection
 singularity, and the 
weighted dual graphs of $V_0$ and $V_t$ are the same.
By \thmref{t:la}, the family $\{V_t | t \in \C\}$ is 
 an equisingular  deformation.

Next assume that $s\ge 2$.
Let $\beta\: Z\times \C \to \C^m\times \C$ be the trivial family of
the weighted blowing up  $Z \to \C^m$ with
 respect to the weight $\w_1=(w_1, \ldots ,w_m)$.
The family $\{V_t|t \in \C\}$ is naturally embedded in
 $\C^m\times\C$. 
Let $W_i=\C^m$ be as above; however we write $x_i$ instead of $z_i$.
Then the cyclic group $\cal C_i$ acts on $W_i\times \C$ as
 in \lemref{l:G}.
Let $V_t^i\subset W_i\times {\{t\}}$ be the strict transform
 of  $V_t$. 
Recall that  $V_t\cap \{x_1=x_2=0\}=\{o\}$ by
 \corref{c:CI} and that the action of the cyclic group $\cal C_i$
 on $V_t^i$  is free outside  the exceptional locus. 
Thus to prove that $V_t$ has an  isolated singularity at the
 origin, it  suffices to show that any component of singular
 loci of $V_t^1$ and  $V_t^2$, 
 intersecting the exceptional set,   is an isolated point.  
Note that the exceptional divisor is singular at  singular
 points of $V_t^i$.   
By \thmref{t:la}, \lemref{l:es-div} and \ref{l:G}, it is
 enough to prove the following three claims (the 
 claims  on $V_t^2$ are proved in the same way).

\begin{clm}\label{clm:1}
  The family of exceptional divisors $E_t \subset V_t^1$ is trivial.
\end{clm}
\begin{clm}\label{clm:2}
 Let $p \in E_0$ be a singular point and let $p_t \in E_t$
 denote the  point $p$        under the identification $E_0=E_t$.
Then each $(V_t^1, p_t)$ is
       an isolated singularity and the weighted dual graph of
 $(V_t^1,x_t)$ is  independent of $t$.
\end{clm}
\begin{clm}\label{clm:3}
 If $q_1, \ldots ,q_k \in V_0^1$ be the fixed points of
 $\cal C_1$-action, 
 then the fixed locus of the family $\{V_t^1|t \in \C\}$ is 
$\bigcup _j\{q_j\}\times \C$. 
If $\cal C_{1,j} \subset \cal C_1$ denotes
 the isotropy group of  $\{q_j\}\times \C$, then
 \lemref{l:G} applies to 
 the family $\{(V_t^1,q_j)|t \in \C\}$ with $\cal C_{1,j}$-action
 for every $q_j$.
\end{clm}

The exceptional divisor $E_t$ is defined by $x_1=0$ in $V_t^1$.
The ideal of  $V_t^1 \subset \C^m= W_1\times {\{t\}}$ is
 generated by  the following functions: 
\begin{align*}
F_{11} & 
=1+a_{11} \x_{d_1-1}^{\alpha_{d_1-1}}+b_{11}
 \bm_{1d_1}+t\bar f_{11}^+, & & \\ 
F_{1j} & 
=\x_j^{\alpha_j}+a_{1j} \x_{d_1-1}^{\alpha_{d_1-1}}+b_{1j}
 \bm_{1d_1} 
 +t\bar f_{1j}^+,   & & 2  \le j \le d_1-2,  \\
F_{ij} & 
=\bm_{ij}+a_{ij}\bm_{id_i-1}+b_{ij}\bar \bm_{id_i}+t\bar f_{ij}^+, 
   & & i  \ge 2, \; 1 \le j \le d_i-2, 
\end{align*}
where 
$$
 \bar f_{ij}^+=f_{ij}^+(x_1^{w_1}, x_1^{w_2}x_2,
 \ldots ,x_1^{w_m}x_m)/x_1^{\word 1(\f_{ij})}, 
$$
and $\bar
 \bm_{id_i}$ is obtained in 
 the same way from $\bm _{id_i}$. 
Recall the condition on the order of $f_{ij}^+$, and that 
$$
 \adeg 1(m_{ij})<\adeg 1(m_{id_i}), \quad i \ge 2, \; 1 \le
 j \le d_i-2. 
$$
Then we see that the ideal of $E_t \subset \C^m$ is
 generated by $x_1$  and the  following polynomials:
\begin{align}
&1+a_{11} \x_{d_1-1}^{\alpha_{d_1-1}}+b_{11}
 \bm_{1d_1},\label{eq:tF1} & & \\ 
&\x_j^{\alpha_j}+a_{1j} \x_{d_1-1}^{\alpha_{d_1-1}}+b_{1j}
 \bm_{1d_1},  & & 2 \le j \le d_1-2, \label{eq:tFj}\\
&\bm_{ij}+a_{ij}\bm_{id_i-1},
  & & i\ge 2, \; 1 \le j \le d_i-2. \label{eq:gij}
\end{align}
These polynomials do not contain $t$; thus  \clmref{clm:1}
 is verified. 
Since the fixed points of the $\cal C_j$-action lie on $E_t$,
 the first assertion of \clmref{clm:3} follows.
By looking  at the action explicitly,
we see that the $\cal C_{1,j}$-action on $W_1\times \C$  is
 unitary around $q_j$.
Thus \clmref{clm:3} follows from \clmref{clm:2}.

It is easy to see that there are $m-d_1$ polynomials in
 \eqref{eq:gij} and any of them 
 contains no variables $x_1, \ldots ,x_{d_1-1}$.
As in the proof of \lemref{l:curveC}, 
we can show that the functions of \eqref{eq:gij} define a
 complete  intersection curve $C' \subset \C^{m-d_1+1}$
 which is smooth  except for the origin. 
Furthermore \eqref{eq:tF1} and \eqref{eq:tFj} define  a
 tower of  cyclic coverings  over $C'$. 
Hence the singularities of $E_t$ are lying above the point
 $(0,\ldots ,0)$  of $C'$.
Therefore at each  singular point of $E_t$, the only $d_1-2$
 variables $x_2,\ldots ,x_{d_1-1}$ 
 are  nonzero,  and others are  zero. 

We fix $t \in \C$. Let $p$ be a singular point of $E_t$.
At $p$,   the Jacobian matrix 
$$
 \partial (F_{11}, \ldots ,F_{1d_1-2})/\partial(x_2, \ldots
 ,x_{d_1-1}) 
$$
 is regular, and  thus  $x_2, \ldots ,x_{d_1-1}$ can be
 expressed by a  convergent power series  with nonzero
 constant terms 
 in $x_1, x_{d_1}, \ldots ,x_m$.  
By substituting them into $F_{ij}$, $i\ge 2$,
 we obtain new defining functions for the germ $(V_t^1, p)$ in
 $(\C^{m-d_1+2}, o)$ as follows:
\begin{equation}\label{eq:neweq}
 F'_{ij}:=\bm_{ij}+a_{ij}\bm_{id_i-1}+b_{ij}'\bar
  \bm_{id_i}'+r_{ij},  
\quad i\ge 2, \quad 1 \le j \le d_i-2,
\end{equation}
where $\bar \bm_{id_i}'$  is a monomial obtained by substituting $1$
 into $x_2, \ldots, x_{d_1-1}$ of $\bar \bm_{id_i}$.
We see  that 
\begin{equation}\label{eq:neq}
 a_{ij_1}b'_{ij_2}\not=a_{ij_2}b'_{ij_1} \; (j_1\neq j_2)
\quad \text{and} \quad  b'_{ij} \neq 0.
\end{equation}

By \eqref{eq:deg}, for every $\Q$-monomial  $m=x(D)$, we have
$$
 \wdeg 1(\bm)=e_1\cdot \adeg 1(m)=e_1\cdot m_{A_1}(D).
$$ 
Suppose that $m_{id_i}=x(D_i)$ for $i\ge 2$.
Since $\adeg 1(m_{id_i-1})=m_{A_1}(\du i)$,  
the exponent of $x_1$ in $\bar \bm_{id_i}'$ is 
\begin{equation}\label{eq:exponent}
  \wdeg 1(\bm_{id_i})-  \wdeg 1(\bm_{id_i-1})=
e_1m_{A_1}(D_i-\du i)>0.
\end{equation}
Since $D_i-\du i \in A_{\Z}$, it follows 
that  $\bar \bm_{id_i}'$ is a monomial of 
$x_1^{e_1}, \x_{d_1}, \ldots ,\x_m$. 
We denote by $h_{i}$ the monomial 
 obtained by replacing  $x_1^{e_1},\x_{d_1}, \ldots ,\x_m$ of 
$\bar \bm_{id_i}'$ with $x_1, x_{d_1}, \cdots ,x_m$, respectively. 
Let $\delta'=(e_1,\delta_{d_1},\ldots ,\delta_m) \in
 \N^{m-d_1+2}$ and let 
$$
f'_{ij}=m_{ij}+a_{ij}m_{id_i-1}+b_{ij}'h_{i}, 
\quad i\ge 2, \quad 1 \le j \le d_i-2.
$$
Then the polynomial
$\bm_{ij}+a_{ij}\bm_{id_i-1}+b_{ij}'\bar \bm_{id_i}'$ is the
 $\delta'$-lifting of $f'_{ij}$.

Let $A_{e}$ be the component of the branch $C_{d_1}$ of the
 node $A_1$,  
which intersects  $A_1$; see \figref{f:diag}.
\begin{figure}[htbp]
  \begin{center}
\setlength{\unitlength}{0.5cm}
 \begin{picture}(14,3)(-7,-2)
\put(0,-1.2){\framebox(6,2){}}
\put(-2,0){\line(1,0){4}}
\put(-2,0){\ten} 
\put(-2.2,-1){$A_1$}
\put(0.8,-1){$A_e$}
\put(2.3,0){ . . .}
\put(1,0){\ten} 
\put(6.5,-0.5){$C_{d_1}$}
\put(-4,-1){\line(2,1){2}}
\put(-6,-2){\framebox(2,1){}}
\put(-8,-2){$C_{d_1-1}$}
\put(-4,1){\line(2,-1){2}}
\put(-6,1){\framebox(2,1){}}
\put(-7.5,1){$C_1$}
\put(-5,-0.5){$\vdots$}
 \end{picture}
\caption{\label{f:diag}}
  \end{center}
\end{figure}
We may assume  that the  $A_e$ is an end of $C_{d_1}$; 
take the blowing up at $A_1 \cap A_e$ if
 necessary.
Let $A'=C_{d_1}$. 
Then $A'\subset M$ can be blown down to a  normal surface 
 singularity. 
We  consider admissible forms concerning $A'$.
We associate the end $A_e$ with the variable $x_1$, while
 keeping the  correspondence between the 
other ends and the variables 
$x_{d_1}, \ldots, x_m$.

\begin{clm}\label{clm:smallgraph}
We have the following.
\begin{enumerate}
 \item $A'$ has $s-1$ nodes $A_2, \cdots, A_s$.
 \item $A'$ satisfies  \condref{c:A},
 and for each $2\le i \le s$,  the set of polynomials
$$
 \{f'_{i1}, \ldots ,f'_{i d_i-2}\}
$$
 is a Neumann-Wahl system concerning $A'$ at  a node $A_i$.
 \item In this new situation, we define a weight $\w_i' \in
       \N^{m-d_1+2}$ 
       ($2\le i \le s$) with respect to $\delta'$  as in the
       sentence before \defref{d:wdeg}. 
Then 
$$ 
\wpdeg i(\bm_{ij}+a_{ij}\bm_{id_i-1}+b_{ij}'\bar \bm_{id_i}') 
<\wpord i(r_{ij}).
$$
\end{enumerate}
\end{clm}

The inductive hypothesis and  \clmref{clm:smallgraph} imply
that $(V_t^1,p)$  is an isolated singularity
and is an equisingular
 deformation of the singularity defined by the $\delta'$-lifting of
the Neumann-Wahl system 
 $\{f'_{ij}\}$; thus \clmref{clm:2} follows.

Now we prove \clmref{clm:smallgraph}.
First, the assertion (1) is obvious.
It is clear that the number of branches of a node $A_i$ ($i\ge 2$)
in $A'$ is the same  as that of $A_i$ in $A$.
For $D=\sum a_jA_j \in A_{\Q}$, we write $\res(D)=
\sum _{A_j\le A'}a_jA_j$.
Let $\t A_k$ denote the dual cycle of $A_k$ on $A'$.
Let $D_{ij}$ ($i \ge 2$) be the monomial
 cycle such that $x(D_{ij})=m_{ij}$, 
and let 
$$
E_{ij}=\res(D_{ij}-\du i).
$$
Since  $(D_{ij}-\du i)\cdot A_e=0$, we have $-E_{ij}\cdot
 A_e=m_{A_1}(D_{ij}-\du i)$; this is the exponent of $x_1$
 in $h_{i}$ 
(see \eqref{eq:exponent}).
By computing the intersection numbers $(E_{ij}+\t A_i)\cdot
 A_l$ for  every $A_l \le A'$,
 we see that $E_{ij}+\t A_i$ is a monomial cycle belonging
 to a branch  of $A_i$ and that
$$
x(E_{id_i}+\t A_i)=h_{i} \quad \text{and} \quad x(E_{ij}+\t
 A_i)=m_{ij} \quad \text{for $1\le j <d_i$}. 
$$
By \eqref{eq:neq}, we have  (2) of \clmref{clm:smallgraph}.
We can generalize the argument above as follows.
Let $m$ be any monomial in $\f_{ij}+tf_{ij}^+$.
Suppose that $m$ is 
 the $\delta$-lifting of a $\Q$-monomial $x(D)$ associated with a
 $\Q$-monomial  cycle $D \in A_{\Q}$.
Let $m'$ be the $\delta'$-lifting of the monomial
 $x(\res(D-\du i)+\t  A_{i})$. 
Then, by the operation which changes
 $\f_{ij}+tf_{ij}^+$ into the function $F'_{ij}$, the monomial 
  $m$  is changed into a function of the form $u m'$,
 where $u \in S$ is a unit. 
We have that  $m_{A_i}(\res(D-\du i)+\t A_{i})=m_{A_i}(D)$
for a node $A_i$ $(i \ge 2)$.
Now it is  easy to see that  (3) holds.
Thus we have proved \clmref{clm:smallgraph}.
\end{proof}

Let $\hat S=\C[[x_1, \ldots ,x_m]]$ denote the  formal power series
ring. 
By  a similar argument as in the proof of \thmref{t:V}, we
obtain the following.

\begin{thm}\label{t:formal}
 Let $\f_{ij}$ be as above. 
Take a formal power series $g_{ij} \in \hat S$ such that
$$
 \wdeg i(\f_{ij})<\word i(g_{ij})\quad \text{ for} \quad
1 \le i\le  s\;, 1 \le j \le d_i-2.
$$ 
Let $J \subset \hat S$ be the ideal generated by all
 $\f_{ij}+g_{ij}$'s. 
Then $\hat S/J$ is a two-dimensional complete intersection ring,
and has only a singularity at the maximal ideal.
\end{thm}

\section{The main results}\label{s:mainresults}

In this section we will prove the following.

\begin{thm}\label{t:main}
  If $\X$ is a rational or minimally elliptic singularity, then its
 universal abelian cover $\Y$ is an equisingular deformation of a 
 Neumann-Wahl complete intersection singularity.
The deformation is defined by the functions of the form
 $f+tf^+$, where $\{f\}$ is a Neumann-Wahl
 system associated with the exceptional set $A$,
 $f^+$ is a  function with  order greater than the degree of
 $f$ and $t$ is the parameter.
\end{thm}

We start without the assumption of the theorem.
We use  the notation of \sref{s:pre}.
Write $\A_b=H^0(-L^{(b)})$
and $\mathcal O_{Y,o}=\bigoplus _{b \in \B}\A_b$.
Recall that for any component $A_i$ of $A$, 
there exist a divisor  $L^i$ and an element $b^i \in \B$ such that 
$\nu (L^i)=\du i$ and $L^i - L^{(b^i)} \in A_{\Z}$.
We consider the
following condition, which depends on  the resolution $\p$.

\begin{cond}\label{c:B}
For each end $A_i \in \E(A)$, 
there exists a section $y_i \in H^0(-L^i)$ such that
 $(y_i)_A=\nu(L^i)$.
\end{cond}

\begin{lem}\label{l:satisfyB}
 If $\X$ is rational,  or if $\X$ is  minimally elliptic and the
 minimally elliptic cycle $E$ is supported on $A$,
then \condref{c:B} is satisfied. 
\end{lem}
\begin{proof}
Let $D$ be any $\pi$-nef divisor on $M$.
If $\X$ is rational, then $D$ is $\pi$-free; see
 \cite[4.17]{chap}.
Assume that $\X$ is  minimally elliptic and
the  minimally elliptic cycle $E$ is supported on $A$.
Since $E$ is 2-connected,   
$H^0(D)$ has no fixed  component on $A$ if $D\cdot A\ge 1$ by
 \cite[4.23, Remark]{chap}. 
Thus the assertion follows.
\end{proof}

Let $\m_Y$ (resp. $\m_X$) denote the maximal ideal of  $\mathcal O_{Y,o}$
(resp. $\mathcal O_{X,o}$).
We may identify $\m_Y$ as $\m_X \bigoplus (\bigoplus _{b
\neq 0}\A_b)$ (cf. \cite[\S 6]{o.uac-rat}).

\begin{lem}\label{l:powerofmax}
 Let $\{Z_k \in A_{\Z}|k \in \N\}$ be a sequence of cycles such that
 $Z_{k+1}>Z_k>0$ for every $k \in \N$.
Then there exists a function $\alpha\:\N \to \N$ such that for 
each $b \in \B$, 
$$
 H^0(-L^{(b)}-Z_k)\subset \m_Y^{\alpha(k)} \quad \text{and}
\quad \lim_{k \to \infty} \alpha (k)= \infty.
$$
\end{lem}
\begin{proof}
We only give an outline.
We can take a positive integer $a$ so that for any $\pi$-nef
 divisor $D$  on $M$ and a cycle 
$Z:=a \sum_{A_i \le A} \du i$,  the natural map
$$
H^0(D-Z)\otimes H^0(-Z) \to H^0(D-2Z)
$$ 
is surjective (cf. \cite[III]{la.simul}).
Let $\beta$ be a nonnegative integer such that $-L^{(b)}-\beta Z$ is
 $\pi$-nef for every $b \in \B$.
Then we obtain $H^0(-L^{(b)}-(\beta+k)Z)\subset \m_Y^k$.
We may assume $Z_1>(\beta+1)Z$.
Now define $\alpha(l)= \max\{k \in \N |(\beta+k)Z \le Z_l\}$.
\end{proof}

\begin{ass}
 From now on, we assume that \condref{c:C} and 
\ref{c:B} are satisfied.
\end{ass}

If $A$ is a chain of curves, then $\X$ is a cyclic quotient
singularity and $\mathcal O_{Y,o}=\C\{y_1,y_2\}$, 
where $y_i$'s are as in \condref{c:B} 
(if $A$ is irreducible, then $y_1,y_2 \in H^0(-L^1)$). 
Let $m=\#\E(A)$. 
Assume that $m\ge 3$.
Then we can define admissible monomials at each node
 by associating each end with a variable as in 
\sref{s:NWS}. 
We define the homomorphism 
$$
\psi\: S=\C\{x_1, \ldots ,x_m\} \to \A=\mathcal O_{Y,o}
$$ 
of $\C$-algebras by $\psi (x_i)=y_i$.
We denote by $\hat{}$ 
the maximal-ideal-adic completions of local rings.
Let $\hat \psi\:\hat S=\C[[x_1, \ldots ,x_m]] \to \hat {\mathcal O}_{Y,o}$ 
be the
induced homomorphism. 
By the definition of the set $\B$,
we may regard  $\B$ as the discriminant
group $G:=\du {\Z}/A_{\Z}$ in the natural way.
Let $S_b \subset S$ (resp. $\hat S_b \subset \hat S$), 
$b \in \B$, denote the set of power series
 represented as the sum of monomials $x(D)$ satisfying 
$D\pmod{A_{\Z}}=b$.
Then we have $S=\bigoplus _{b \in \B}S_b$ and the $\psi$ becomes a
homomorphism of $\B$-graded (or $G$-graded) algebras.
The same holds for $\hat S$ and $\hat{\psi}$.

Let $A_1, \ldots ,A_s$ be all of the nodes of $A$, 
and $d_i$  the number of branches of a node $A_i$.
Let $\M_i=\{m_{i1}, \ldots ,m_{id_i}\}$ denote a complete system of
admissible monomials at a node $A_i$. 
Let $\C\M_i \subset S$ denote the $\C$-linear subspace spanned by
the monomials of $\M_i$.

\begin{lem}\label{l:constructCSAF}
 For any node $A_i$, let $\mu_i$ denote the composite of
 homomorphisms 
$$
\C\M_{i} \overset{\psi}\longrightarrow 
H^0( \mathcal O_M(-L^{i})) \to H^0(\mathcal O_{A_{i}}(-L^{i})).
$$
Then $\mu_i$ is surjective.
We have  $h^0(\mathcal O_{A_{i}}(-L^{i}))=2$ and $\dim \Ker \mu_i=d_i-2$.
Let  $\F_i=\{f_{i1}, \ldots ,f_{id_i-2}\}$ be a basis of 
$\Ker \mu_i$. 
Then $\F_i$ is a Neumann-Wahl system at the node $A_i$.
\end{lem}
\begin{proof}
Since $L^i\cdot A_i=-1$, we have $h^0(\mathcal O_{A_{i}}(-L^{i}))=2$.
Suppose $m_{ij}=x(D_{ij})$ for a monomial cycle $D_{ij}$ belonging 
to a  branch $C_{ij}$ of $A_i$.
Since $(D_{ij}-\du {i})\cdot A_{i}=1$,
it follows from \lemref{l:product} that 
 $\mu_i(x(D_{ij}))$ has a zero of order one at  $A_{i} \cap C_{ij}$.
Thus $\mu_i(x(D_{id_i-1}))$ and $\mu_i(x(D_{id_i}))$ generate
 $H^0(\mathcal O_{A_{i}}(-L^{i}))$.
Therefore we obtain a complete system of admissible forms
 expressed by 
a $((p-2) \times p)$-matrix as in \remref{r:normal}, which
 is a basis of 
 $\Ker \mu_i$.
\end{proof}

Let us recall that polynomials of $\F_i$ are quasihomogeneous with
respect to the  $A_i$-weight. 

\begin{lem}\label{l:higherterms}
 Let $A_i$ be a node and $h\in H^0(-L^i)$. 
Then there exists $\bar h \in \hat S_{b^i}$
 such that $\hat \psi(\bar h)=h$ in $\hat {\mathcal O}_{Y,o}$.
Suppose that $h=\psi (\bar h_0)$ for an
 admissible form $\bar h_0 \in \Ker \mu _i$.
Then we can take the $\bar h$ so that 
$\adeg i (\bar h_0)<\aord i(\bar h)$.
If in addition $\psi $ is surjective, 
such $\bar h$ can be taken from $S_{b^i}$.
\end{lem}
\begin{proof}
We use the notation of \lemref{l:constructCSAF}.
Suppose that $h \neq 0$.
 Let $F_0$ be the divisor such that $\nu (F_0)=(h)_A$ and
 $F_0-L^i \in  A_{\Z}$.   
Let $c_0=-F_0\cdot A_i$. Then $c_0\ge 0$.
By \condref{c:C} and the proof of \lemref{l:existmonomials},
 there exists a cycle  $F'_0\ge F_0$ such that
 $m_{A_i}(F'_0-F_0)=0$, $(F'_0-F_0)\cdot A_i=0$ and
 $F'_0\cdot A_j=0$ if  $j\neq i$ and $A_j$ is not an end.
Let $D'_{ij}=D_{ij}-\du i$. 
Then an arbitrary cycle of the form $E_{\mathbf
 a}:=F'_0+\sum a_{j}D'_{ij}$ is a    monomial cycle, 
where $\mathbf a=(a_1, \ldots ,a_m) \in \Z^m$ with 
$\sum a_j=c_0$ and $a_j \ge 0$. 
It is easy to see that the $\mu_i(x(E_{\mathbf a}))$'s span  
$H^0(\mathcal O_{A_i}(-F_0))$. 
Thus we have a quasihomogeneous polynomial $\bar h_1$,
which is  a linear form of $x(E_{\mathbf a})$'s, with respect to
 $A_i$-weight such that $h-\psi(\bar h_1) \in H^0(-F_0-A_i)$.
Let $F_1$ be the divisor such that $\nu (F_1)=(h-\psi(\bar h_1))_A$ 
and $F_1-L^i \in A_{\Z}$.  
Then it follows from the  argument above that there exists
 a quasihomogeneous polynomial $\bar h_2$ 
 such that 
$$
h-\psi(\bar h_1)-\psi(\bar h_2) \in H^0(-F_1-A_i).
$$ 
Thus we obtain a sequence $\{\bar h_k|k \in \N\}$
of quasihomogeneous  polynomials  and a sequence 
$\{F_k | k  \in \N\}$ 
of  divisors  satisfying the following: for all $k \in \N$,
\begin{enumerate}
 \item $h-\psi(\bar h_1+\cdots +\bar h_k) \in H^0(-F_k)$,
 \item $F_{k+1}>F_k$,
 \item $\adeg i(\bar h_k)<\adeg i(\bar h_{k+1})$.
\end{enumerate}
By \lemref{l:powerofmax}
there exists a function $\alpha\: \N \to \N$ such
 that $H^0(-F_k) \subset \m_Y^{\alpha(k)}$ and  
 that $\lim_{k\to \infty}\alpha(k) = \infty$.
Now put $\bar h=\sum \bar h_i \in \hat S_{b^i}$.
Then $\hat \psi(\bar h)=h$.
Suppose that $h=\psi(\bar h_0)$ with an admissible form
$\bar h_0 \in \Ker  \mu_i$.
Since $\psi(\bar{h_0}) \in H^0(-L^i-A_i)$, 
we have 
$$
 \adeg i(\bar h_0)<m_{A_i}(F_0)=\adeg  i(\bar h_{1}).
$$
If $\psi$ is surjective, then the maps $(x_1, \cdots ,x_m)^k
 \to \m_Y^k$ and $S_{b} \to \A_b$   are surjective 
for every $k \in \N$ and $b \in \B$. 
Therefore, for sufficiently large $k$, there exists $\bar h'
 \in S_{b^i}$ such that $\aord i(\bar h')>\adeg i(\bar h_0)$
 and $h=\psi(\bar h_1+\cdots +\bar h_k+\bar h')$.
\end{proof}

\begin{prop}\label{p:surj}
 The homomorphism $\psi\: S \to \mathcal O_{Y,o}$ is surjective.
\end{prop}
\begin{proof}
We fix a node $A_i$.
By \cite[Proposition 5.1]{nw-CIuac}, the group
$G=\du {\Z}/A_{\Z}$ is generated by
 $\{\du j|A_j \in \E(A)\}$.
Hence for each $b \in \B$, 
there exists a monomial $m_b \in S$ such that
 $\A_b\cdot m_b \subset  H^0(-L^i)$.
By \lemref{l:higherterms}, 
we have $\A_b \subset \hat\psi(\hat S_b)\cdot
 m_b^{-1}$. Therefore
$$
\hat \psi (\hat S) \subset \hat {\mathcal O}_{Y,o} 
\subset \sum_{b \in \B} \hat \psi (\hat S) \cdot m_b^{-1}.
$$
Then it follows that $\hat S /\Ker \hat \psi$ is a
 two-dimensional domain. 
By \lemref{l:higherterms}, for any $f_{kl} \in \bigcup \F_j$,
there exists $\tilde f_{kl}  \in \hat S$ such  that
 $\LF_{A_k}(\tilde f_{kl})=f_{kl}$  
and $\tilde f_{kl} \in \Ker \hat \psi$.
Let $\tilde I \subset \hat S$ denote the ideal generated by 
all $\tilde f_{kl}$'s. 
Then it follows from \thmref{t:formal} 
that $\hat S/\tilde I$ is a two-dimensional normal domain. 
Since $\tilde I \subset \Ker \hat{\psi}$,
we obtain that  $\hat S/\tilde I\cong \hat \psi (\hat S)$.
Since $\hat S$ is Noetherian, $ \hat {\mathcal O}_{Y,o} $ is finitely
 generated  $\hat \psi(\hat S)$-module.
By the normality of $\hat \psi (\hat S)$,
we obtain $\hat \psi (\hat S)=\hat {\mathcal O}_{Y,o}$.
This implies that $\psi$ is surjective.
\end{proof}

It follows from \lemref{l:higherterms} and \proref{p:surj}
that for each 
$f_{ij} \in \F_i$,  there exists $f_{ij}^+ \in S_{b^i}$ 
such that $f_{ij}+f_{ij}^+
\in \Ker \psi$ and $\adeg i( f_{ij})<\aord i(f_{ij}^+)$.
For each $t \in \C$, we denote by $I_t \subset S$ the ideal
generated by the functions
$$
f_{ij}+tf_{ij}^+, \quad 1\le i \le s, 1 \le j \le d_{i}-2.
$$
As in the proof of \proref{p:surj}, by \thmref{t:V}, 
 we see that $\Ker \psi=I_1$.

\begin{cor}
 $\mathcal O_{Y,o} \cong S/I_1$.
\end{cor}

Again by \thmref{t:V}, we obtain the following.

\begin{thm}\label{t:main-general}
 If \condref{c:C} and \ref{c:B} are satisfied, then the universal
 abelian cover  $\Y$ of $\X$
 is an equisingular deformation of a 
 Neumann-Wahl complete intersection singularity.
The deformation is defined by the functions 
$$
f_{ij}+T f_{ij}^+ \in S[T], \quad  1\le i \le s, 1 \le j \le
 d_{i}-2, 
$$
where $S[T]$ is the polynomial ring over $S$.
\end{thm}

Now \thmref{t:main} follows from \thmref{t:main-general},
\lemref{l:existmonomials} and
\ref{l:satisfyB}. 

\begin{cor}\label{c:QGor}
  If \condref{c:C} and \ref{c:B} are satisfied, then  $\X$ is 
 $\Q$-Gorenstein.
If in addition the link of $\X$ is a  homology sphere,
then $\X$ is an equisingular deformation of a 
 Neumann-Wahl complete intersection singularity.
\end{cor}

\begin{rem}
 There is ambiguity in the choice of the complete systems of
 admissible monomials $\M_i$. 
 The proof of \lemref{l:higherterms} shows that 
 if $m$ and $m'$ belong to the same branch of a node $A_i$,
 then there exists $h \in S_{b^i}$ such that $x(m)-x(m')-h \in \Ker
 \psi$ and $\aord i(h)>\adeg i(x(m))$. 
 Therefore, whether $\Y$ is a Neumann-Wahl complete
 intersection depends not only on the choice of the sections
 $\{y_i\}$, but also on the choice of the complete
systems of admissible monomials. 
\end{rem}

In the rest of this section, we describe the action of the
Galois group of the universal abelian covering.
Recall that the Galois group  $H_1(\Sigma, \Z)$ of the universal
abelian covering $\Y \to \X$  is isomorphic
to  the discriminant group $G=\du {\Z}/A_{\Z}$.

For any $\Q$-cycle $D=\sum a_i\du i=\sum b_iA_i$, 
we have $D\cdot A_j=-a_j$ and $D\cdot \du j=-b_j$.
Thus we obtain the following

\begin{lem}\label{l:du}
Let $D \in A_{\Q}$. Then
\begin{enumerate}
 \item $D \in \du {\Z}$ if and only if $\{D\cdot A_i|A_i \le
       A\} \subset \Z$, 
 \item $D \in A_{\Z}$ if and only if $\{D\cdot \du i|A_i \le
       A\}  \subset \Z$. 
\end{enumerate}
 
\end{lem}
For any $D \in \du {\Z}$, we denote by $(D)$ the class of
$D$ modulo $A_{\Z}$. 
Then the action of $G$ on $S=\C\{x_1,\ldots  ,x_m\}$ is
expressed as follows (cf. \cite[\S 5]{nw-CIuac}). 
For $(D) \in G$ and a monomial cycle $F$, 
\begin{equation}\label{eq:action}
  (D)\cdot x(F):=\exp (2\pi\sqrt{-1}D\cdot F)x (F).
\end{equation}
By \lemref{l:du} (1), this is well-defined.
Let $(Y_t,o)$ denote the singularity defined by 
the ideal $I_t \subset S$.  
We have seen that $\{Y_t|t \in \C\}$ is an equisingular family.
Since  $I_t$ is generated by homogeneous elements
of $G$-graded algebra $S$, the group $G$ naturally acts on $S/ I_t
=\mathcal O_{Y_t,o}$.
It follows from \lemref{l:du} (2) that $\mathcal
O_{X,o}=(\mathcal O _{Y,o})^{G}$.
By \corref{c:CI} and \cite[Proposition 5.2]{nw-CIuac},
 the action is free on $Y_t\setminus \{o\}$ (see also \cite[Theorem
 7.2 (2)]{nw-CIuac}).  
Therefore  $X_t:=Y_t/G$ is a normal  singularity.
The linear action of $G$ on $\C^m$ determined by
\eqref{eq:action} is unitary.
By \lemref{l:G} and the uniqueness of the universal abelian
covering, we obtain the following

\begin{thm}\label{p:esquotient}
The family $\{X_t|t \in \C\}$ is an equisingular
 deformation, and each $Y_t \to X_t$ is the universal abelian
covering.
\end{thm}


\begin{thebibliography}{10}

\bibitem{bri.simul}
E.~Brieskorn, \emph{Singular elements of semi-simple algebraic groups}, Actes
  du Congr\`es International des Math\'ematiciens (Nice, 1970), Tome 2,
  Gauthier-Villars, Paris, 1971, pp.~279--284.

\bibitem{is.simul}
S.~Ishii, \emph{Simultaneous canonical models of deformations of isolated
  singularities}, Algebraic geometry and analytic geometry (Tokyo, 1990),
  ICM-90 Satell. Conf. Proc., Springer, Tokyo, 1991, pp.~81--100.

\bibitem{la.rat}
H.~Laufer, \emph{On rational singularities}, Amer. J. Math. \textbf{94} (1972),
  597--608.

\bibitem{la.me}
H.~Laufer,  \emph{On minimally elliptic singularities}, Amer. J. Math. \textbf{99}
  (1977), 1257--1295.

\bibitem{la.simul}
H.~Laufer,  \emph{Weak simultaneous resolution for deformations of {G}orenstein
  surface singularities}, Singularities, Part 2 (Arcata, Calif., 1981), Proc.
  Sympos. Pure Math., vol.~40, Amer. Math. Soc., Providence, R.I., 1983,
  pp.~1--29.

\bibitem{neumann.plumbing}
W.~D. Neumann, \emph{A calculus for plumbing applied to the topology of complex
  surface singularities and degenerating complex curves}, Trans. Amer. Math.
  Soc. \textbf{268} (1981), 299--344.

\bibitem{neumann.abel}
W.~D. Neumann, \emph{Abelian covers of quasihomogeneous surface singularities},
  Singularities, Part 2 (Arcata, Calif., 1981), Proc. Sympos. Pure Math.,
  vol.~40, Amer. Math. Soc., Providence, RI, 1983, pp.~233--243.

\bibitem{nw-uac}
W.~D. Neumann and J.~Wahl, \emph{Universal abelian covers of surface
  singularities}, Trends in singularities, Trends Math., Birkh\"auser, Basel,
  2002, pp.~181--190.

\bibitem{nw-HSL}
W.~D. Neumann and J.~Wahl, \emph{Complex surface singularities with integral homology sphere
  links}, 2003, preprint {(math.AG/0301165)}.

\bibitem{nw-qcusp}
W.~D. Neumann and J.~Wahl, \emph{Universal abelian covers of quotient-cusps}, Math. Ann.
  \textbf{326} (2003), 75--93.

\bibitem{nw-CIuac}
W.~D. Neumann and J.~Wahl, \emph{Complete intersection singularities of splice type as universal
  abelian covers}, 2004, preprint {(math.AG/0407287)}.

\bibitem{o.uac-rat}
T.~Okuma, \emph{Universal abelian covers of rational surface singularities}, J.
  London Math. Soc. (2) \textbf{70} (2004), 307--324.

\bibitem{chap}
M.~Reid, \emph{Chapters on algebraic surfaces}, Complex algebraic geometry,
  IAS/Park City Math. Ser., vol.~3, Amer. Math. Soc., Providence, RI, 1997,
  pp.~3--159.

\bibitem{tki-w}
M.~Tomari and K.{-i}. Watanabe, \emph{Filtered rings, filtered blowing-ups and
  normal two-dimensional singularities with ``star-shaped'' resolution}, Publ.
  Res. Inst. Math. Sci. \textbf{25} (1989), 681--740.

\bibitem{wahl.es}
J.~Wahl, \emph{Equisingular deformations of normal surface singularities, {I}},
  Ann. of Math. \textbf{104} (1976), 325--365.

\bibitem{wahl.defqh}
J.~Wahl, \emph{Deformations of quasi-homogeneous surface singularities}, Math.
  Ann. \textbf{280} (1988), 105--128.

\end{thebibliography}

\end{document}